\documentclass[12pt]{article}
\usepackage{amsmath,amsthm,amsfonts}
\theoremstyle{plain}
\newtheorem{theor}{Theorem}
\theoremstyle{plain}
\newtheorem{lemma}{Lemma}
\theoremstyle{plain}
\newtheorem{cor}{Corollary}
\theoremstyle{remark}
\newtheorem{rem}{Remark}
\normalbaselines
\begin{document}
\vbox {\vspace{6mm}}
\begin{center}
{\bf\large Disjointly Strictly Singular Inclusions\\
of Symmetric Spaces}\\[3mm]
{\bf S.V.Astashkin}\\
{\it Samara State University, 443011 Samara 11, Russia}\\[5mm]
\end{center}
\begin{abstract}
\noindent
In this paper, the disjoint strict singularity of inclusions of
symmetric spaces of functions on an interval is considered. A condition for
the presence of a "gap" between spaces sufficient for the inclusion
of one of these spaces into the other to be disjointly strictly singular
is found. The condition is stated in terms of fundamental functions of spaces
and is exact in a certain sense. In parallel, necessary and sufficient
conditions for an inclusion of Lorentz spaces to be disjointly strictly
singular (and similar conditions for Marcinkiewicz spaces) are obtained
and certain other assertions are proved.
\\
\\
\end{abstract}

\begin{center}
{\bf Introduction}\\[1mm]
\end{center}

Recall that a bounded linear operator
$T$
from a Banach space
$X$
into a Banach space
$Y$
is called {\it strictly singular} (or a {\it Cato operator}) if
$X$
does not contain an infinite-dimensional subspace
$Z$
such that the restriction of
$T$
to
$Z$
is an isomorphism.

In recent decades the class of strictly singular operators has been
extensively studied (see the references in, e.g., the monograph [1]).
One of the historically first results important for our purposes
is the Grothendick theorem on strict singularity of the identity inclusion
operator from
$L_{\infty}(\Omega,\mu)$
into
$L_p(\Omega,\mu),$
where
$1\leq p<\infty$
and
$\mu$
is a probability measure on
$\Omega$
(see [2] or [3, Theorem 5.2]). However, as a rule, the identity inclusion
operator from one symmetric space into another (the definition
is given below) is not strictly singular because
of the existence of "through" subspaces (such as the subspace generated
by the Rademacher functions [4]). In part because of this, the close
notion of disjointly strictly singular operator was introduced in 1989 [5].

A bounded linear operator
$T$
from a Banach lattice
$X$
into a Banach space
$Y$
is called {\it disjointly strictly singular} (or {\it has the DSS property})
if there exists  no sequence of nonzero disjoint vectors
$\{x_n\}_{n=1}^{\infty}$
in
$X$
such that the restriction of
$T$
to their closed linear hull
$[x_n]$
is an isomorphism.

Clearly, any strictly singular operator is a DSS operator. A simple
example shows that the converse is not true. For instanse,
the identity inclusion operator
$I: L_p[0,1]\to{L_q[0,1]}$
$(1\leq q<p\leq{\infty})$
has the DSS property, because the closed linear hull in
$L_r$
of disjoint functions
$x_n\in{L_r[0,1]}$
is isomorphic to
$l_r$
$(1\leq r\leq{\infty}):$
$\;[x_n]_r\approx {l_r}.$
However, if
$p<\infty,$
then Khintchine's inequality [6] implies that
$[r_n]_p\approx{[r_n]_q}\approx{l_2}$
($r_n$
are Rademacher functions), and therefore
$I$
is not strictly singular. At the same time, it is easy to show that if
$X$
has a Schauder basis of disjoint vectors, then the class of DSS operators
on
$X$
coincides with the class of strictly singular operators [7].

The notion of DSS operator proved important in studies of the
geometric properties of function spaces. For example, the existence of
operators without the DSS property makes it possible to construct
complemented subspaces that admit "nonstandard" projections onto them [5,7].

The goal of this paper is to study the following question: when does
the identity inclusion operator (throughout, we denote it by
$I$)
from one symmetric space in another have the DSS property? The conditions
are stated in terms of fundamental functions of these spaces.

If
$z=z(t)$
is measurable on
$[0,1]$ with respect to the Lebesgue measure
$\mu$,
then we call the function
 $n_{z}\,(\tau)=\mu\{t:\,|z(t)|>\tau\}$
$(\tau>0)$
the {\it distribution function} of $z.$
Two functions
$x(t)$
and
$y(t)$
are called {\it equimeasurable} if
$n_x(\tau)=n_y(\tau)$
for
$\tau>0.$

Recall that a Banach space $E$ of measurable functions on
$[0,1]$
is called a symmetric space (briefly is an SS) if the following conditions
hold:

(1) if
$y\in{E}$
and
$|x(t)|\le{|y(t)|},$
then
$x\in{E}$
and
$||x||\le{||y||};$

(2) if
$y\in{E}$
and
functions
$x(t)$
and
$y(t)$
are equimeasurable,
then $x\in{E}$
and
$||x||=||y||.$

The fundamental function of an SS
$E$
is defined by
$f_E(t)=||\chi_{(0,t)}||_E,$
where, as usual,
$\chi_U(t)=1\,(t\in U),$
$\chi_U(t)=0\,(t\not\in U).$
The function
$f_E(t)$
is quasiconcave on
$(0,1]$
[8, p.137], i.e., it is nonnegative, increasing, and
$f_E(t)/t$
decreases. As is known (see, e.g., [8, p.70]), such a function is
equivalent to its least concave majorant. Throughout,
$G$
denotes the class of all positive increasing functions concave on
$(0,1].$

An important example of an SS is an Orlicz space. Let
$N(t)$
be an increasing convex function on
$[0,\infty)$
such that
$N(0)=0$
and
$N(\infty)=\infty.$
The Orlicz space
$L_N$
consists of all functions
$x=x(t)$
measurable on
$[0,1]$
and such that
$$
\int_T{N\left(\frac{|x(t)|}{u}\right)\,d{\mu}}\,<\,\infty$$
for some
$u>0;$
the norm of this space is
$$
||x||\,=\,\inf\left\{u>0:\;\int_T{N\left(\frac{|x(t)|}{u}\right)\,d{\mu}}\,\leq 1\right\}.$$

Direct calculation shows that the fundamental function of the space
$L_N$
is
$f_N(t)=1/N^{-1}(1/t)$
($N^{-1}(u)$
is the inverse of
$N(u)$)
[9].
\vskip 0.2cm

In [5], the following disjoint strict singularity theorem for inclusions
of
$L_N$
into
$L_M$
is proved.

\noindent {\bf Theorem.} {\it If
$L_N\subset{L_M},$
then the following conditions are equivalent:

(1) the inclusion
$I: L_N\to{L_M}$
is a DSS operator;

(2) for any
$n=1,2,..$
and
${\cal K}>0,$
there exist
$1\leq x_1<x_2<...<x_n$
and
$c_1>0,..,c_n>0$
such that
$$
\sum_{i=1}^n\,c_iN(tx_i)\,\geq{\,{\cal K}\sum_{i=1}^n\,c_iM(tx_i)}\;\;\;\;\mbox{for}\;\;t\geq 1.$$}
\vskip 0.3cm

Let us show that condition (2) follows from the relation
$$
\lim_{t\to +0}\,\frac{f_M(t)}{f_N(t)}\,=\,0,$$
where
$f_N$
and
$f_M$
are the fundamental functions of the respective Orlicz spaces.
Indeed, this relation implies that
$N^{-1}(t)\leq{hM^{-1}(t)}$
for an arbitrary positive
$h\leq 1$
and
$t\geq t_0$
and, since
$N(t)$
is convex,
$M(t)\leq{N(ht)}\leq{hN(t)},$
if
$t\geq{M^{-1}(t_0)}.$
Therefore,
$$
\lim_{t\to{\infty}}\,\frac{M(t)}{N(t)}\,=\,0,$$
and condition (2) holds.

Quite naturally, this observation leads us to the following general problem.

\def\ab{I:E\to F}
\def\bc{E\subset F}
\def\cd{{\cal M}_f(t)}
\def\de{\varphi\in G}
\def\ef{\psi\in G}
\def\fg{\Lambda(\varphi)}
\def\gh{\Lambda(\psi)}
\def\hj{\varepsilon}
\def\ji{\frac{\psi(t)}{\varphi(t)}}
\def\ik{M(\theta)}
\def\kl{M(\tilde{\varphi})}
\def\lm{M(\tilde{\psi})}
\def\mn{||w_m||_2}
\def\no{\sum_{k=n_m}^{n_{m+1}-1}}
\def\op{\frac{dt}{\psi(t)}}
\def\pq{\log_2^{1/2}{4/t}}
\def\qr{\max_{0\leq m\leq r}}
\def\rs{\max_{n_m\leq k<n_{m+1}}}
\def\st{\biggl\|\sum_{m=0}^{\infty}\,a_mv_m\biggr\|}
\def\tz{\int_0^1\,\frac{d\psi(s)}{\varphi(s)}}
\def\zx{\sum_{k=0}^{\infty}\,\frac{\psi(2^{-k})-\psi(2^{-k-1})}{\varphi(2^{-k})}}
\def\xy{\sum_{k=0}^{\infty}\,\frac{a_k}{\sqrt{S_k}\varphi(2^{-k})}}
\def\ya{M(\tilde{\rho})\subset{\Lambda(\psi)}}

Suppose that functions
$\de$
and
$\ef$
satisfy the condition

$\;(A)\;\;\lim_{t\to +0}\psi(t)/\varphi(t)\,=\,0$\\
and are, respectively, the fundamental functions of symmetric spaces
$E$
and
$F$
such that
$E\subset F.$
Does the identity inclusion operator
$I:\,E\to F$
have the DSS property?

In what follows, we show that the answer to this question for "classical"
symmetric spaces (such as the Lorentz and Marcinkiewicz spaces) as well
as for Orlicz spaces is positive. Moreover, it is so for an inclusion of
a Lorentz space into an arbitrary SS (and, vice versa, of an arbitrary
SS into a Marcinkiewicz space).

However, in the general case, this is not true: this paper contains an
example of two symmetric spaces
$E$
and
$F$
such that
$\bc$
and their fundamental functions satisfy condition (A), but 
$\ab$
is not a DSS operator.

At the same time, it is possible to state a condition on fundamental
functions stronger than (A) under which the answer to the stated question
is positive for all symmetric spaces. First, recall the definition of the
dilation function.

For a positive function
$f$
on
$(0,1],$
its {\it dilation function}
$\cd$
is defined as
$$
\cd\,=\,\sup\left\{\frac{f(st)}{f(s)}:\;0<s\leq\min\left(1,\frac{1}{t}\right)\right\}\;\;\;(t>0).$$
Since
$\cd$
is semimultiplicative, there exist numbers
$$
\gamma_f=\,\lim_{t\to 0}\frac{\ln{\cd}}{\ln t}\;\;\;\mbox{and}\;\;\;\delta_f=\,\lim_{t\to \infty}\frac{\ln{\cd}}{\ln t},$$
which are called, respectively, the {\it lower} and {\it upper dilation indices}
of function
$f$.
If
$\de$,
then we have
$0\leq \gamma_{\varphi}\leq \delta_{\varphi}\leq 1$
[8, p.76].

Let us introduce one more condition on the functions
$\varphi$
and
$\psi$
from the class $G:$

$\;(B)\;\;\gamma_{\psi/\varphi}\,>\,0.$\\
The definition of lower dilation index readily implies that condition
(A) follows from (B). The converse, of course, is not true: it suffices
to take for
$\varphi$
and
$\psi$
functions differing by a logarithmic factor (see also the proof of
Theorem 3).

We shall show that, if (B) holds, then operator
$I:\,E\to F$
will be a DSS operator for arbitrary symmetric spaces
$E$
and
$F,$
$E\subset F,$
with fundamental functions
$\varphi$
and
$\psi,$
respectively. This result generalizes and simultaneously refines
a similar theorem for the Orlicz spaces proved in [7]. In parallel, we
shall show that condition (A) is necessary and sufficient for the
identity inclusion operator from one Lorentz space into another to have
the DSS property (and similar assertion for Marcinkiewicz spaces).
These results also supplement the theorem for Orlicz spaces proved in [5]
and cited above.
\vskip 0.6cm

\begin{center}
{\bf ${\cal x}\,1.$  The inclusions} $\fg\subset F$ {\bf and} $E\subset \ik$
\end{center}
\vskip 0.2cm

For
$\de$,
the Lorentz space
$\fg$
consists of all functions
$x=x(s)$
measurable on
$[0,1]$
and such that
$$
||x||_{\fg}\,=\,\int_0^1\,x^*(s)\,d\varphi(s)\,<\,\infty\,,$$
where
$x^*(s)$
is a decreasing left-continuous rearrangement of the function
$|x(s)|$
[8, p.83]. Clearly, the fundamental function of the Lorentz space
$\Lambda(\varphi)$
is
$f_{\fg}(t)=\,\varphi(t).$

\begin{theor}
Suppose that the functions
$\de$
and
$\ef$
satisfy (A) and
$F$
is an SS on
$[0,1]$
with fundamental function
$\psi(t).$
Then
$\fg\subset F$
and the identity inclusion
$I:\,\fg\to F$
is a DSS operator.
\end{theor}
\begin{proof}
Condition (A) and the continuity
of the concave functions
$\varphi$
and
$\psi$
at
$t>0$
imply that
$\psi(t)\leq{C_1\varphi(t)}$
for
$t\in{[0,1]}.$
By the definition of the norm of a Lorentz space, we have
$\fg\subset \gh.$
The Lorentz space with a given fundamental function is minimal among
the symmetric spaces with the same fundamental function [8, p.160];
therefore,
$\fg\subset F.$

Suppose that
$I:\,\fg\to F$
has not the DSS property. Then there exists a sequence of nonzero
disjoint functions
$x_n\geq 0$
such that
$$
||x_n||_{\fg}\,\leq{C_2||x_n||_F}\;\;\mbox{for}\;\;n=1,2,..\eqno{(1)}$$

By condition (A), for any
$0<\hj<1,$
there exists an
$h>0$
such that
$$
\psi(t)\leq{\hj\varphi(t)}\eqno{(2)}$$
for all positive
$t<h.$
Choose
$N$
so that, for
$n\geq N,$
$\mu(g_n)<h,$
where
$g_n\,=\,\{\,t\in{[0,1]}\,:\,x_n(t)\ne 0\}.$

The subset of finite-valued functions is dense in any Lorentz space on
$[0,1]$
[8, p.149].
Therefore, for each
$n\geq N,$
there exists a function
$$
y_n(t)\,=\,\sum_{k=1}^{m_n}a_k^n\chi_{e_k^n},\;\;\;\mbox{where}\;\;a_k^n\geq 0,\;e_1^n\supset{e_2^n}\supset{...}\supset{e_{m_n}^n},\;\;\;\mbox{and}\;\;\;\mu(e_1^n)<h,$$
for which
$$
\max(||x_n-y_n||_{\fg},||x_n-y_n||_{\gh})\,<\,\hj \min(||x_n||_{\fg},||x_n||_{\gh}).\eqno{(3)}$$
Hence
$||x_n||_{\gh}-||y_n||_{\gh}\,\leq{\hj ||x_n||_{\gh}},$
and [8, p.160] and (2) imply that
\begin{multline*}
||x_n||_F\leq{||x_n||_{\gh}}\leq{1/(1-\hj)\,||y_n||_{\gh}}\,=\,1/(1-\hj)\,\sum_{k=1}^{m_n}a_k^n\psi(\mu(e_k^n))\leq\\
\leq{\hj /(1-\hj)\,\sum_{k=1}^{m_n}a_k^n\varphi(\mu(e_k^n))}\,=\,\hj /(1-\hj)||y_n||_{\fg}.
\end{multline*}
In addition, by (3),
$$
||y_n||_{\fg}\,\leq{\,(1+\hj)||x_n||_{\fg}}.$$
Thus,
$$
||x_n||_F\,\leq{\,\frac{\hj(1+\hj)}{1-\hj}||x_n||_{\fg}}.$$

This inequality with an
$\hj>0$
satisfying
$$
\frac{1-\hj}{\hj(1+\hj)}\,>\,C_2$$
contradicts (1).
\end{proof}
\begin{cor}
Suppose that
$\de$,
$\ef$,
and
$\psi(t)\leq{C_1\varphi(t)}$
for
$t\in{(0,1]}.$
The following conditions are equivalent:

(1) (A) holds;

(2) the inclusion
$I:\,\fg\to{\gh}$
is a DSS operator;

(3) there exist no sequence of nonzero disjoint functions
$\{x_n\}$
and no constant
$C_2>0$
such that
$$
||x_n||_{\fg}\,\leq{\,C_2||x_n||_{\gh}}\;\;\;\mbox{for}\;\;n=1,2,...\eqno{(4)}$$
\end{cor}
\begin{proof}
The implications
$(1)\to{(2)}$
and
$(2)\to{(3)}$
follow from Theorem 1 and the definition of a DSS operator, respectively.

Suppose that condition (A) does not hold, i.e., that
$$
\limsup_{t\to 0}\frac{\psi(t)}{\varphi(t)}\,>\,0.$$
Then there exist a sequence
$\{t_k\}\subset{(0,1]}$
and a constant
$C_2>0$
such that we have
$\sum_{n=1}^{\infty}t_n\leq 1$
and
$\varphi(t_n)\leq{C_2\psi(t_n)}$
for
$n=1,2,...$
The functions
$x_n=\chi_{e_n},$
where
$e_n\subset{[0,1]}$
are disjoint and
$\mu(e_n)=t_n,$
satisfy (4). Therefore, (3) implies (1); this completes the proof of
Corollary 1.
\end{proof}

Let
$\theta$
be a function from
$G.$
The Marcinkiewicz space
$M(\theta)$
consists of all functions
$x=x(s)$
measurable on
$[0,1]$
and such that
$$
||x||_{\ik}\,=\,\sup_{0<t\leq 1}\frac{1}{\theta(t)}\,\int_0^t{x^*(s)ds}\,<\,\infty.$$
The fundamental function of the space
$\ik$
is
$f_{\ik}(t)=\tilde{\theta}(t)=t/\theta(t).$

\begin{theor}
Let the functions
$\de$
and
$\ef$
satisfy condition (A). If
$E$
is an SS on
$[0,1]$
with fundamental function
$\varphi(t),$
then
$E\subset{\lm}$
and the identity inclusion
$I:\,E\to{\lm}$
is a DSS operator.
\end{theor}
\begin{proof}
Since any Marcinkiewicz space is maximal among all symmetric spaces with
the same fundamental function [8, p.162], we have
$E\subset{\lm}.$

By condition (A),
$\psi(t)\leq{C_1\varphi(t)}$
for
$t\in{(0,1]};$
therefore, by the definition of a Marcinkiewicz space,
$\kl\subset{\lm}.$
Hence,
$E\subset{\lm}.$

Suppose that
$I:\,E\to{\lm}$
is not a DSS operator. Then, in particular, there exists a sequence
of disjoint functions
$x_n$
such that
$$
||x_n||_{\lm}=1\;\;\mbox{and}\;\;||x_n||_{\kl}\leq C_2\;\;\;\mbox{for}\;\;n=1,2,...\eqno{(5)}$$
Choose a
$t_k\in{(0,1]}$
for which
$$
\int_0^{t_k}x_k^*(s)ds\,\geq{\,\frac{1}{2}\tilde{\psi}(t_k)}\;\;\;\;(k=1,2,..).$$
Since the functions
$x_k$
are disjoint, we can assume that
$t_k\to 0.$
Therefore,
$$
||x_k||_{\kl}\;\geq\;{\frac{1}{\tilde{\varphi}(t_k)}\,\int_0^{t_k}x_k^*(s)ds}\;\geq\;{\frac{\varphi(t_k)}{2\psi(t_k)}}.$$
By (A),
$||x_k||_{\kl}\to{\infty}$ as
$k\to{\infty},$
which contradicts condition (5).

This completes the proof of Theorem 2.
\end{proof}
\begin{cor}
Suppose that
$\de$,
$\ef$,
and
$\psi(t)\leq{C_1\varphi(t)}$
for
$t\in{(0,1]}.$
The following conditions are equivalent:

(1) (A) holds;

(2) the inclusion
$I:\,\kl\to{\lm}$
is a DSS operator;

(3) there exist no sequence of nonzero disjoint functions
$x_n$
and no constant
$C_2>0$
such that
for some $C_2>0$
$$
||x_n||_{\kl}\,\leq{\,C_2||x_n||_{\lm}}\;\;\;\mbox{for}\;\;n=1,2,...$$
\end{cor}

The proof of Corollary 2 is similar to the proof of Corollary 1.

\begin{cor}
For an arbitrary SS
$E\ne{L_1}$
on
$[0,1],$
the inclusion
$I:\,E\to{L_1}$
is a DSS operator.
\end{cor}
\begin{proof}
First,
$L_1=M(1)$
and an arbitrary SS
$E$
is embedded in
$L_1$
[8, p.124]. If
$f_E(t)=\varphi(t),$
then the function
$t/\varphi(t)$
increases, because
$\varphi$
is concave. Therefore, condition (A) is violated if and only if
$\varphi(t)\approx t$
(i.e., if and only if
$C_1t\leq{\varphi(t)}\leq{C_2t}$
for some
$C_1>0$
and
$C_2>0$),
and
$E=L_1.$
It remains to apply Theorem 2.
\end{proof}
\begin{rem}
The assertion of Corollary 3 was proved in [10] in a different way.
\end{rem}
\begin{rem}
Arguing as in the proof of Corollary 3 and applying Theorem 1, we can
readily show that the inclusion
$I:\,L_{\infty}\to{E}$
is a DSS operator for any SS
$E\ne{L_{\infty}}.$
Moreover, it is shown in [11] that this operator is even strictly singular.
This generalizes the Grothendieck theorem mentioned in the introduction.
\end{rem}

In the next section we show that, generally, condition (A) is not
sufficient for the inclusion of an SS with fundamental function
$\varphi$
into an SS with fundamental function
$\psi$
to have the DSS property.
\vskip 0.6cm

\begin{center}
{\bf ${\cal x}\,2.$\,\,An example of symmetric spaces
$E$
and
$F$
such that
$E\subset F$\\
and their fundamental functions satisfy condition (A),\\
but
$I:\,E\to F$
is not a DSS operator}
\end{center}
\begin{theor}
There exist two symmetric spaces
$E$
and
$F$
on
$[0,1]$
with fundamental functions
$\varphi$
and
$\psi,$
respectively, such that
$\bc,$
$\varphi$
and
$\psi$
satisfy condition (A), and the operator
$\ab$
has not the DSS property.
\end{theor}
\begin{proof}
Let the SS
$E$
be the Marcinkiewicz space
$\lm$
with
$\tilde{\psi}(t)=t/{\psi(t)},$
where
$$
\psi(t)=\,t^{1/2}\,\log_2^{1/2}\frac{4}{t},\;\;\;0<t\leq 1.$$
It is readily verified that
$\psi$
is an increasing concave function on
$[0,1]$
and
$\gamma_{\psi}=\delta_{\psi}=1/2.$
Therefore, by [8, p.156],
$$
||x||_{\lm}\,\approx{\,\sup_{0<t\leq 1}\{x^*(t)\psi(t)\}}.\eqno{(6)}$$

Let us define the space
$F$.
We put
$b_k=\,(k+2)^{-1/2}2^{k/2}$
and
$z_k(t)=\,b_k\chi_{(0,2^{-k}]}(t)$
and define a sequence of numbers
$n_0=1<n_1<n_2<...<n_m<...$
by setting
$$
n_{m+1}=\,\max\{n=1,2,..:\,\sum_{k=n_m}^{n-1}\frac{1}{k+2}\leq 1\}\eqno{(7)}$$
and a sequence of functions
$w_m=w_m(t)$
by setting
$$
w_m(t)=\,\max_{n_m\leq k<n_{m+1}}z_k(t),\;\;\;m=0,1,..$$
Since
$\{b_k\}$
increases, (7) implies that the norms of
$w_m$
in
$L_2$
satisfy the inequalities:
$$
||w_m||_2^2\,\geq{\,\sum_{k=n_m}^{n_{m+1}-1}b_k^22^{-k-1}}\,=\,
1/2\no\frac{1}{k+2}\,\geq\,\frac{1}{4}$$
and
$$
{\mn}^2\,\leq{\,{\no}b_k^22^{-k}}\,=\,\no\frac{1}{k+2}\,\leq\,1.$$
Therefore,
$$
\frac{1}{2}\leq\mn\leq 1.\eqno{(8)}$$

Consider
$\chi_b=\,b^{-1/2}\chi_{(0,b)},$
$\bar{w}_m=\,w_m/{||w_m||_2},$
and
$V\,=\,\{\chi_b\}_{0<b\leq 1}\bigcup{\{\bar{w}_m\}_{m=0}^{\infty}}.$
Let
$F$
be the set of all functions
$x=x(t)$
measurable on
$[0,1]$
and satisfying
$$
||x||=\,\sup\left\{\int_0^1x^*(t)v(t)dt:\,v\in V\right\}\,<\,\infty.$$
Then
$F$
is an SS on
$[0,1]$
as the intersection of Lorentz spaces determined by the functions
$\int_0^tv(s)ds$
with
$v\in V.$
In addition, the definition of
$F$
implies that
$||x||_{M(t^{1/2})}\leq{||x||_F}\leq{||x||_2}.$
Therefore the fundamental functions
$f_E(t)=\psi(t)$
and
$f_F(t)=t^{1/2}$
of the spaces
$E$
and
$F$
satisfy condition (A).

Let us prove that
$$
\lm\subset F.\eqno{(9)}$$
By (6), it is suffices to show that
$1/\psi\in F.$
Indeed,
$$
\int_0^1\chi_b(t)\op\,=\,2b^{-1/2}\int_0^b\frac{d(t^{1/2})}{\pq}\,\leq 2\;\;\mbox{for}\;\;0<b\leq 1,$$
$$
\int_0^1w_m(t)\op\,\leq{\,{\no}b_k\int_0^{2^{-k}}\op}\,=\,2{\no}b_k\int_0^{2^{-k}}\frac{d(t^{1/2})}{\pq}\leq$$
$$
\leq{\,2{\no}\frac{1}{k+2}}\,\leq\,2.$$
Therefore (8) and the definition of
$F$
imply that
$||1/\psi||_F\leq 4,$
which proves (9).

Next, we put
$D_m=(2^{-n_{m+1}},2^{-n_m}]$
and
$$
v_m(t)\,=\,w_m(t)\chi_{D_m}(t)\,=\,{\no}b_k\chi_{(2^{-k-1},2^{-k}]}(t)\;\;\mbox{for}\;\;m=0,1,...$$
The functions
$v_m$
are disjoint. Let us show that the norms of
$E$
and
$F$
are equivalent on their linear hull.

Suppose that
$$
v(t)=\,\sum_{m=0}^ra_mv_m(t).$$
Without loss of generality, we can assume that
$a_m\geq 0.$
Consider
$w(t)=\,{\qr}a_mw_m(t).$
The function
$w(t)$
monotonically decreases on
$(0,1],$
and
$v(t)\leq w(t).$
Therefore, by (6),
$$
||v||_E\,\leq\,||w||_E\,\leq{\,C\qr\left\{a_m{\rs}b_k\psi(2^{-k})\right\}}.$$
Since
$b_k\psi(2^{-k})=1$
for
$k=0,1,2,..,$
we obtain
$$
||v||_E\,\leq{\,C{\qr}a_m}.\eqno{(10)}$$

Now, let us estimate
$||v||_F$
from below. By (7), for any
$m=0,1,..,r$
we have
\begin{multline*}
\int_0^1v_m^*(t)w_m(t)dt\,\geq{\,\int_0^1(v_m^*(t))^2dt}\,=\,\int_0^1v_m^2(t)dt\,=\\
={\no}b_k^22^{-k-1}\,=\,\frac{1}{2}\no\frac{1}{k+2}\,\geq\,\frac{1}{4}.
\end{multline*}
Hence (8) implies that
$||v_m||_{F}\geq{1/4}.$
Therefore, by (10),
$$
||v||_F\,\geq{\,\qr\{a_m||v_m||_F\}}\,\geq{\,\frac{1}{4}{\qr}a_m}\,\geq{\,\frac{||v||_E}{4C}}.$$
This together with (9) means that the norms of the spaces
$E$
and
$F$
are equivalent on the linear hull of the set of functions
$v_m$
$(m=0,1,..).$
Hence there exists a
$B>0$
such that, for an arbitrary
$a_m,$
$$
B^{-1}{\st}_F\,\leq{\,{\st}_E}\,\leq{\,B{\st}_F}.$$
In other words, the identity inclusion operator
$\ab$
has not the DSS property.

This completes the proof of Theorem 3.
\end{proof}
\begin{rem}
Theorem 3 shows that relation (A) does not guarantee the presence of a
"gap" between spaces sufficient for the corresponding identity inclusion
operator to have the DSS property. On the other hand, simple examples
show that, even for symmetric spaces with the same fundamental function,
this operator may have this property.

For instance, the inclusion of the Lorentz space
$\Lambda(t^{1/p})$
into the space
$L_p,$
where
$1<p<\infty,$
has the DSS property. Indeed, it is easy to show that any sequence of
normalized disjoint functions in the Lorentz space contains a subsequence
equivalent to the standard basis in
$l_1$.
At the same time, any such sequence in
$L_p$
is equivalent to the standard basis in
$l_p.$
\end{rem}

In the last section we show that, unlike (A), condition (B) is sufficient for
the inclusion operator
$I:\,E\to F$
of arbitrary symmetric spaces
$E$
and
$F$
with fundamental functions
$\varphi$
and
$\psi,$
respectively, to be a DSS operator.
\vskip 0.6cm

\begin{center}
{\bf ${\cal x}\,3.$\,\,Suffiiency of condition (B) for the operator}\\
$I:\,E\to F$
{\bf to have the DSS property}
\end{center}
\begin{theor}
Suppose that functions
$\de$
and
$\ef$
satisfy condition (A),
$\delta_{\varphi}<1,$
and we have
$\kl\subset{\gh}.$
Then the operator
$I:\,\kl\to{\gh}$
has the DSS property.
\end{theor}

First, we prove the following auxiliary assertion.
\begin{lemma}
Under the assumptions of Theorem 4, there exists a function
$\rho\in G$ such that
$$
(1)\;\lim_{t\to 0}\frac{\rho(t)}{\varphi(t)}\,=\,0;$$
$$
(2)\;\ya.$$
\end{lemma}
\begin{proof}
Since
$\delta_{\varphi}<1,$
[8, p.156] implies that
$$
||x||_{\kl}\,\approx{\,\sup_{0<t\leq 1}\,\{\varphi(t)x^*(t)\}}.$$
Therefore the relation
$\kl\subset{\gh}$
is equivalent to
$$
\tz\,<\,\infty.\eqno{(11)}$$
Since the function
$\varphi$
is concave, we have
$$
\tz\,\approx{\,\zx};$$
hence (11) is equivalent to the condition
$$
\zx\,<\,\infty.\eqno{(12)}$$
Put
$$
a_k=\,\psi(2^{-k})-\psi(2^{-k-1})\;\;\;\mbox{and}\;\;\;S_n=\,\sum_{k=n}^{\infty}\frac{a_k}{\varphi(2^{-k})}.$$
Then
$S_n\to 0,$
and [12, Chap. 3, Ex. 12] and (12) imply that
$$
\xy\,<\,\infty\eqno{(13)}$$
By the definition of upper dilation index, there exist
$u>0$
and
$C>0$
such that
$\delta_{\varphi}+u<1$
and
$$
{\cal M}_{\varphi}(t)\,\leq{\,C\,t^{\delta_{\varphi}+u/2}}\eqno{(14)}$$
for all
$t\geq 1.$
Consider the sequence of numbers
$$
g_0=\,S_0,\;\;g_k=\,\max(S_k,2^{-u}g_{k-1})\;\;\mbox{for}\;\;k=1,2,...\eqno{(15)}$$
Suppose that
$g=g(t)$
is a function linear on the intervals
$[2^{-k-1},2^{-k}],$
$g(2^{-k})=g_k$
for
$k=0,1,..,$
and
$h(t)=\,\sqrt{g(t)}\varphi(t).$

Since
$\{S_k\}$
decreases,
$\{g_k\}$
also decreases; therefore, the functions
$g(t)$
and
$h(t)$
increase on
$(0,1].$
It follows from (15) that, for
$j\geq 0$
and
$2^{-k-1}<t\leq{2^{-k}},$
$$
\frac{g(2^jt)}{g(t)}\leq{\frac{g(2^{j-k})}{g(2^{-k-1})}}=\frac{g_{k-j}}{g_{k+1}}=\frac{g_{k+1-(j+1)}}{g_{k+1}}\leq{2^{(j+1)u}}.$$
Hence, by (14),
$$
{\cal M}_h(2^j)\,=\,\sup_{0<t\leq{2^{-j}}}\,\frac{\varphi(2^jt)\sqrt{g(2^jt)}}{\varphi(t)\sqrt{g(t)}}\,\leq{\,C2^{u/2}2^{(\delta_{\varphi}+u)j}}\;\;\mbox{for}\;\;j=0,1,...$$
This implies that
$\delta_h<1,$
because
$\delta_{\varphi}+u<1.$
Therefore, according to [8, p.78], the function
$h(t)$
is equivalent to its least concave majorant; we denote this majorant by
$\rho(t)$.

The function
$\rho$
belongs to
$G,$
and
$\rho(2^{-k})\approx{h(2^{-k})}=\sqrt{g_k}\varphi(2^{-k})$
for
$k=0,1,2,..;$
in addition, since
$g_k\geq{S_k},$
it follows from (13) that
$$
\sum_{k=0}^{\infty}\,\frac{a_k}{\rho(2^{-k})}\,\leq{\,C_1\xy}\,<\,\infty.$$
This equivalent to
$$
\ya;$$
the equivalence is proved in the same way as for the function
$\varphi$
(see (11) and (12)).

Finally, relation (15) gives
$$
\lim_{t\to 0}\frac{\rho(t)}{\varphi(t)}\,\leq{\,C_1\lim_{t\to 0}\sqrt{g(t)}}\,=\,\lim_{k\to{\infty}}\sqrt{g_k}\,=\,0.$$

This completes the proof of Lemma 1.
\end{proof}

\begin{proof}[Proof of Theorem 4.]
By Lemma 1, there exists a function
$\rho\in G$
such that
$$
\kl\subset{\ya}\;\;\;\mbox{and}\;\;\;\lim_{t\to 0}\frac{\rho(t)}{\varphi(t)}\,=\,0.$$
According to Corollary 2, the operator
$I:\,\kl\to{M(\tilde{\rho})}$
has the DSS property; all the more, it has this property when regarded
as an operator from
$\kl$
into
$\gh$.

This proves Theorem 4.
\end{proof}

Theorem 4 makes it possible to prove the sufficiency of condition (B)
for the identity inclusion operator from an SS
$E$
into an SS
$F$
with fundamental functions
$\varphi$
and
$\psi,$
respectively, to have the DSS property.

\begin{theor}
Suppose that functions
$\de$
and
$\ef$
satisfy condition (B) and
$E$
and
$F$
are symmetric spaces with fundamental functions
$\varphi$
and
$\psi,$
respectively. Then
$\bc$
and
$\ab$
is a DSS operator.
\end{theor}
\begin{proof}
Let us verify that the assumptions of Theorem 4 hold, i.e., that
$\delta_{\varphi}<1$
and that
$$
\kl\subset{\gh}.\eqno{(16)}$$

First, condition (B) implies the existence of
$u>0$
and
$C>0$
such that
$$
\frac{\psi(ts)\varphi(s)}{\psi(s)\varphi(ts)}\,\leq{\,Ct^u}\eqno{(17)}$$
whenever
$0<t\leq 1$
and
$0<s\leq 1.$
This and the concavity of
$\psi$
give
$$
{\cal M}_{\varphi}\left(\frac{1}{t}\right)\,\leq{\,C{\cal M}_{\psi}\left(\frac{1}{t}\right)t^u}\,=\,Ct^{u-1};$$
therefore,
$\delta_{\varphi}\leq{1-u}<1.$

Next, since the function
$x^*(t)$
decreases, we have
$$
x^*(t)\leq{\frac{||x||_{\kl}}{\varphi(t)}}\;\;\;\mbox{if}\;\;x\in{\kl}\;\;\mbox{and}\;\;0<t\leq 1.$$
Therefore, to prove (16), it suffices to verify that
$1/\varphi\in{\gh}.$

By (17),
$\psi(t)/\varphi(t)\leq{C_1t^u}.$
Hence
$\psi(t)\to 0$
as
$t\to 0,$
and
$$
||1/\varphi||_{\gh}\,=\,\int_0^1\frac{\psi'(t)}{\varphi(t)}dt\,\leq{\,\int_0^1\frac{\psi(t)}{\varphi(t)}\frac{dt}{t}}\,\leq{\,\frac{C_1}{u}}\,<\,\infty.$$

This proves (16). The above-mentioned extremality of the Lorentz and
Marcinkiewicz spaces in the class of symmetric spaces with the same fundamental
function [8] implies that
$E\subset{\kl}\subset{\gh}\subset F.$
Therefore,
$\ab$
has the DSS property, because it has this property as an operator from
$\kl$
into
$\gh$
by Theorem 4.

This completes the proof of Theorem 5.
\end{proof}

\begin{rem}
Theorem 5 was proved in [7] for Orlicz spaces under a condition on
functions
$\varphi$
and
$\psi$
somewhat more restrictive than (B), namely, the inequality
$\delta_{\varphi}<\gamma_{\psi}.$
\end{rem}

The author wishes to express his gratitude to S. Ya. Novikov for
useful discussiones.
\newpage
\pagestyle{empty}

\end{document}